# EXPLICIT REPRESENTATION OF FINITE PREDICTOR COEFFICIENTS AND ITS APPLICATIONS

By Akihiko Inoue and Yukio Kasahara[1]

*Hokkaido University*

We consider the finite-past predictor coefficients of stationary time series, and establish an explicit representation for them, in terms of the MA and AR coefficients. The proof is based on the alternate applications of projection operators associated with the infinite past and the infinite future. Applying the result to long memory processes, we give the rate of convergence of the finite predictor coefficients and prove an inequality of Baxter-type.

**1. Introduction.** Let $\{X_k\} = \{X_k : k \in \mathbf{Z}\}$ be a real, zero-mean, weakly stationary process defined on a probability space $(\Omega, \mathcal{F}, P)$, which we shall simply call a *stationary process*. We denote by $H$ the real Hilbert space spanned by $\{X_k : k \in \mathbf{Z}\}$ in $L^2(\Omega, \mathcal{F}, P)$. The norm of $H$ is given by $\|Y\| := E[Y^2]^{1/2}$. For $n \in \mathbf{N}$, we denote by $H_{[-n,-1]}$ and $H_{(-\infty,-1]}$ the subspaces of $H$ spanned by $\{X_{-n}, \ldots, X_{-1}\}$ and $\{X_k : k \leq -1\}$, respectively. We write $P_{[-n,-1]}$ and $P_{(-\infty,-1]}$ for the orthogonal projection operators of $H$ onto $H_{[-n,-1]}$ and $H_{(-\infty,-1]}$, respectively. The projection $P_{[-n,-1]}X_0$ (resp., $P_{(-\infty,-1]}X_0$) stands for the best linear predictor of the future value $X_0$ based on the finite past $\{X_{-n}, \ldots, X_{-1}\}$ (resp. the infinite past $\{X_k : k \leq -1\}$), and its mean square prediction error is given by $\sigma_n^2 := \|X_0 - P_{[-n,-1]}X_0\|^2$ (resp. $\sigma^2 := \|X_0 - P_{(-\infty,-1]}X_0\|^2$).

For nondeterministic $\{X_k\}$ (see Section 2.1), the finite predictor coefficients $\phi_{n,j}$ are the uniquely determined ones in

$$(1.1) \qquad P_{[-n,-1]}X_0 = \sum_{j=1}^{n} \phi_{n,j} X_{-j}.$$

Received November 2004; revised June 2005.

[1] Supported in part by the Grant-in-Aid for JSPS Fellows, The Ministry of Education, Culture, Sports, Science and Technology, Japan.

As is well known, we can calculate the numerical values of $\phi_{n,1},\ldots,\phi_{n,n}$, as well as the mean square prediction error $\sigma_n^2$, from the values $\gamma(0),\ldots,\gamma(n)$ of the autocovariance function of $\{X_k\}$, using recursive algorithms such as the Durbin–Levinson algorithm (see, e.g., Section 5.2 in [6]). The recursive methods are of great practical importance in time series analysis. However, they are not necessarily effective in problems of a theoretical character, in particular, those related to the asymptotic behavior as $n \to \infty$.

A classical problem of this type is the rate of convergence of $\sigma_n^2 - \sigma^2 \downarrow 0$ as $n \to \infty$. See, for example, [8], where references to earlier work—by Grenander and Rosenblatt, Grenander and Szegö, Baxter, Ibragimov and many others—are given. The arguments in these references are closely related to the theory of orthogonal polynomials as described in [10, 26, 27].

A new approach to a related problem was introduced by Inoue [15]. For the partial autocorrelation coefficients $\alpha(n) = \phi_{n,n}$ of a stationary process $\{X_k\}$ with short or long memory, the asymptotic behavior of $|\alpha(n)|$ as $n \to \infty$ was obtained using a representation of the mean square prediction error $\sigma_n^2$ in terms of the MA (moving-average) coefficients $c_k$ and the AR (autoregressive) coefficients $a_k$ (see Section 2.2 for the definitions of $c_k$ and $a_k$). By the same approach, but with extra complication, similar results on $|\alpha(n)|$ were obtained in [17, 18] for the fractional ARIMA (autoregressive integrated moving-average) processes. The fractional ARIMA model is an important parametric model including a class of long memory processes. It was independently introduced by Granger and Joyeux [9] and Hosking [13] (see Example 2.6). The advantage of such an approach, that is, that via $c_k$ and $a_k$, has become more apparent in [16] where a representation of the partial autocorrelation function $\alpha(\cdot)$ itself, in terms of $c_k$ and $a_k$, was derived. The representation enabled us to study the behavior of $\alpha(\cdot)$ more directly, and thereby to improve results in several ways. In particular, the asymptotic behavior of $\alpha(n)$ as $n \to \infty$, rather than that of $|\alpha(n)|$, was obtained.

In this paper our main interest is in the finite predictor coefficients $\phi_{n,j}$, which are among the most basic quantities in the prediction theory for $\{X_k\}$. After we establish an explicit representation of the type above for $\phi_{n,j}$, that is, that in terms of the MA coefficients $c_k$ and the AR coefficients $a_k$, we provide two applications of the representation to long memory processes.

For $n \in \mathbf{N}$, we write $H_{[-n,\infty)}$ for the subspace of $H$ spanned by $\{X_k : k \geq -n\}$ and $P_{[-n,\infty)}$ for the orthogonal projection operator of $H$ onto $H_{[-n,\infty)}$. To prove the representation of $\phi_{n,j}$, we use an approximation scheme based on the alternate applications of the projections $P_{(-\infty,-1]}$ and $P_{[-n,\infty)}$. In so doing, the following equalities play a key role:

$$(1.2) \qquad H_{(-\infty,-1]} \cap H_{[-n,\infty)} = H_{[-n,-1]}, \qquad n = 1, 2, \ldots$$

(see Theorem 2.2). For example, it is known that (1.2) holds if $\{X_n\}$ is purely nondeterministic and has spectral density $\Delta(\cdot)$ such that

$$\text{(1.3)} \qquad \int_{-\pi}^{\pi} \frac{1}{\Delta(\lambda)} \, d\lambda < \infty$$

(Theorem 3.1 in [15]). We discuss the equivalence between (1.2) and *complete nondeterminism* (see Theorem 2.3 and Remark 2).

When Wiener's prediction formula (1.4) below is available, we thus obtain a representation of $\phi_{n,j}$ in terms of $c_k$ and $a_k$ (Theorem 2.5). However, in applications it is essential that $\phi_{n,j}$ be expressed in terms of absolutely convergent series made up of $c_k$ and $a_k$. We derive such an expression (Theorem 2.9) under additional conditions on $c_k$ and $a_k$, that is, (A1) or (A2) in Section 2.3. The condition (A1) corresponds to short memory processes, and (A2) to long memory processes.

The first application of the representation of $\phi_{n,j}$ concerns the rate of convergence of $\phi_{n,j}$ toward its limit as $n \to \infty$. Under suitable conditions, $\phi_{n,j}$ converges to the infinite predictor coefficient $\phi_j$ in

$$\text{(1.4)} \qquad P_{(-\infty,-1]}X_0 = \sum_{j=1}^{\infty} \phi_j X_{-j}.$$

The rate at which $\phi_{n,j}$ converges to $\phi_j$ is a fundamental problem in prediction theory and time series analysis. A textbook treatment of this problem can be found in [22], Section 7.6. Using the representation of $\phi_{n,j}$, we show that

$$\phi_{n,j} - \phi_j = \sum_{k=2}^{\infty} g_k(n,j),$$

where $g_k(n,j)$ is a function of $\{c_k\}$ and $\{a_k\}$ [see (2.28)], and we examine the convergence rate for a long memory process whose autocovariance function $\gamma(\cdot)$ is regularly varying at infinity with index $-p$ for some $p \in (0,1)$. It is shown that $\lim_{n\to\infty} n\{\phi_{n,j} - \phi_j\}$ exists, and the limit is calculated exactly in terms of $p$ and $\{\phi_k\}$ (Theorem 3.3). It is interesting that the rate of convergence does not depend on $p$.

The second application of the representation of $\phi_{n,j}$ is related to the additional error $\|P_{[-n,-1]}X_0 - \sum_{j=1}^{n} \phi_j X_{-j}\|$ that arises when we use the infinite predictor coefficients $\phi_j$ instead of the finite ones $\phi_{n,j}$. There exists a known inequality that deals with this problem, and is commonly referred to as *Baxter's inequality* (see [1]; see also [3, 7] and Section 7.6.2 in [22]). It takes the form

$$\text{(1.5)} \qquad \sum_{j=1}^{n} |\phi_{n,j} - \phi_j| \leq M \sum_{k=n+1}^{\infty} |\phi_k|,$$

with finite positive constant $M$. The original inequality (1.5) of Baxter was an assertion for short memory processes. By simple arguments based on the representation of $\phi_{n,j}$, we prove (1.5) for long memory processes, including the fractional ARIMA processes (Theorem 4.1).

In Section 2 we prove the representation of the finite predictor coefficients $\phi_{n,j}$. In Section 3 we apply it to show the rate of convergence of $\phi_{n,j}$ for long memory processes. In Section 4 we apply the representation to prove an inequality of Baxter-type for long memory processes.

**2. Finite predictor coefficients.** Let $\{X_n\} = \{X_n : n \in \mathbf{Z}\}$ be a stationary process; as stated in Section 1, this means that $\{X_n\}$ is a real, zero-mean, weakly stationary process defined on a probability space $(\Omega, \mathcal{F}, P)$. The *autocovariance function* $\gamma(\cdot)$ of $\{X_n\}$ is defined by

$$\gamma(n) := E[X_n X_0], \qquad n \in \mathbf{Z}.$$

As we also stated in Section 1, we denote by $H$ the closed real linear hull of $\{X_k : k \in \mathbf{Z}\}$ with respect to the norm $\|Y\| := E[Y^2]^{1/2}$. Then $H$ is a real Hilbert space with inner product $(Y, Z) := E[YZ]$. For $n, m \in \mathbf{Z}$ with $n \leq m$, we write $H_{(-\infty,n]}$, $H_{[n,\infty)}$, $H_{[n,m]}$ and $H_{\{n\}}$ for the closed subspaces of $H$ spanned by $\{X_k : -\infty < k \leq n\}$, $\{X_k : n \leq k < \infty\}$, $\{X_k : n \leq k \leq m\}$ and $X_n$, respectively. Notice that $H_{\{n\}} = H_{[n,n]}$. For an interval $I$, we write $P_I$ for the orthogonal projection operator of $H$ onto $H_I$.

A stationary process $\{X_n\}$ is said to be *purely nondeterministic* (PND) if

$$\bigcap_{n=-\infty}^{\infty} H_{(-\infty,n]} = \{0\}.$$

If there exists an even, nonnegative and integrable function $\Delta(\cdot)$ on $[-\pi, \pi]$ such that

$$\gamma(n) = \int_{-\pi}^{\pi} e^{in\lambda} \Delta(\lambda) \, d\lambda, \qquad n \in \mathbf{Z},$$

then $\Delta(\cdot)$ is called a *spectral density* of $\{X_n\}$. As is well known, $\{X_n\}$ is PND if and only if it has a positive spectral density such that $\int_{-\pi}^{\pi} |\log \Delta(\lambda)| \, d\lambda < \infty$ (see, e.g., Chapter II in [23]).

2.1. *Convergence of an approximation scheme.* Let $Y \in H$. If $\{X_k\}$ is *nondeterministic*, that is, $X_0 \notin H_{(-\infty,-1]}$, then $X_{-n}, \ldots, X_{-1}$ are linearly independent, whence we can express the predictor $P_{[-n,-1]}Y$ uniquely in the form

$$(2.1) \qquad P_{[-n,-1]}Y = \sum_{j=1}^{n} \phi_{n,j}(Y) X_{-j}.$$

In this section we prove the convergence of an approximation scheme for computing the real coefficients $\phi_{n,j}(Y)$.

For $n, k \in \mathbf{N}$, we define the orthogonal projection operator $P_n^k$ by

$$(2.2) \qquad P_n^k := \begin{cases} P_{(-\infty,-1]}, & k = 1, 3, 5, \ldots, \\ P_{[-n,\infty)}, & k = 2, 4, 6, \ldots. \end{cases}$$

It should be noticed that $\{P_n^k : k = 1, 2, \ldots\}$ is merely an alternating sequence of projection operators, first to the subspace $H_{(-\infty,-1]}$, then to $H_{[-n,\infty)}$, and so on.

LEMMA 2.1. *Assume that $\{X_n\}$ is nondeterministic. Let $Y$ be an arbitrary element of $H$. Then, for $n, k \in \mathbf{N}$, there exist unique real coefficients $\phi_{n,1}^k(Y), \ldots, \phi_{n,n}^k(Y)$, as well as $Z_n^k \in H_{(-\infty,-n-1]}$ for $k$ odd and $Z_n^k \in H_{[0,\infty)}$ for $k$ even, such that*

$$P_n^k P_n^{k-1} \cdots P_n^1 Y = \sum_{j=1}^n \phi_{n,j}^k(Y) X_{-j} + Z_n^k.$$

PROOF. We assume that $k$ is odd. From Lemma 6.1 in [22] (Regression Lemma), it follows that

$$(2.3) \qquad H_{(-\infty,-1]} = H_{(-\infty,-n-1]} + H_{[-n,-1]} \qquad \text{(direct sum)}$$

(see the proof of Theorem 6.3 in [22]). Since $X_{-n}, \ldots, X_{-1}$ are linearly independent and $P_n^k P_n^{k-1} \cdots P_n^1 Y \in H_{(-\infty,-1]}$, the lemma for $k$ odd follows. The case in which $k$ is even is proved in a similar fashion. □

It is natural to ask if $\phi_{n,j}^k(Y)$ converges to $\phi_{n,j}(Y)$ as $k \to \infty$.

THEOREM 2.2. *We assume that*

$$(2.4) \qquad \{X_n\} \text{ is nondeterministic and satisfies } (1.2).$$

*Then we have*

$$(2.5) \qquad \phi_{n,j}(Y) = \lim_{k \to \infty} \phi_{n,j}^k(Y), \qquad Y \in H, \ n \in \mathbf{N}, \ j = 1, \ldots, n.$$

*In particular, (2.5) holds if*

$$(2.6) \qquad \{X_n\} \text{ is purely nondeterministic and satisfies } (1.3).$$

PROOF. The condition (1.2) and von Neumann's alternating projection theorem (see, e.g., Theorem 9.20 in [22]) yield

$$(2.7) \qquad \operatorname*{s-lim}_{k \to \infty} P_n^k P_n^{k-1} \cdots P_n^1 = P_{[-n,-1]}, \qquad n = 1, 2, \ldots.$$

We put

(2.8) $$\varepsilon_k := X_k - P_{(-\infty,k-1]}X_k \qquad (k \in \mathbf{Z}).$$

Then, from Lemma 2.1 we see that

$$(P_n^{2k+1}P_n^{2k}\cdots P_n^1 Y, \varepsilon_{-1}) = \phi_{n,1}^{2k+1}(Y) \cdot \|\varepsilon_{-1}\|^2.$$

By (2.7), the left-hand side tends to $(P_{[-n,-1]}Y, \varepsilon_{-1})$ as $k \to \infty$. Thus, $a_{n,1} := \lim_{k\to\infty} \phi_{n,1}^{2k+1}(Y)$ exists. In the same way, letting $k \to \infty$ in

$$(P_n^{2k+1}P_n^{2k}\cdots P_n^1 Y, \varepsilon_{-2}) = \phi_{n,2}^{2k+1}(Y) \cdot \|\varepsilon_{-2}\|^2 + \phi_{n,1}^{2k+1}(Y) \cdot (X_{-1}, \varepsilon_{-2}),$$

we find the existence of $a_{n,2} := \lim_{k\to\infty} \phi_{n,2}^{2k+1}(Y)$. Repeating this argument, we see that $a_{n,j} := \lim_{k\to\infty} \phi_{n,j}^{2k+1}(Y)$ exists for all $j = 1, \ldots, n$. Hence, $Z_n := \lim_{k\to\infty} Z_n^{2k+1}$ also exists in $H$, and we have

$$Z_n = P_{[-n,-1]}Y - \sum_{j=1}^n a_{n,j}X_{-j}.$$

Since the right-hand side is in $H_{[-n,-1]}$, so is $Z_n$. Moreover, $Z_n \in H_{(-\infty,-n-1]}$ since, for every $k \geq 1$, $Z_n^{2k+1}$ belongs to the closed subspace $H_{(-\infty,-n-1]}$. Combining, $Z_n \in H_{[-n,-1]} \cap H_{(-\infty,-n-1]}$. However, by (2.3) this implies $Z_n = 0$. Thus, $P_{[-n,-1]}Y = \sum_{j=1}^n a_{n,j}X_{-j}$. By uniqueness, we obtain $\phi_{n,j}(Y) = a_{n,j} = \lim_{k\to\infty} \phi_{n,j}^{2k+1}(Y)$. Similarly, we have $\phi_{n,j}(Y) = \lim_{k\to\infty} \phi_{n,j}^{2k}(Y)$. Thus, (2.5) follows. Finally, by Theorem 3.1 in [15], (2.6) implies (2.4), whence (2.5). □

REMARK 1. A stationary process $\{X_n\}$ is said to be *minimal* if $X_0$ does not belong to the closed linear span of $\{X_k : k \in \mathbf{Z}, \ k \neq 0\}$ in $H$. By Theorem 24 in [20], (2.6) is equivalent to saying that $\{X_n\}$ is purely nondeterministic and minimal. The condition (2.6) is also equivalent to another property called *pure minimality* (see [21, 24] and Theorem 8.10 in [22]). The condition (2.6) holds in most interesting examples, and we can easily check it.

Since the assumption (2.4) is a key to our arguments, we are interested in its characterization. The next theorem gives such a result.

THEOREM 2.3. *The condition* (2.4) *is equivalent to*

(2.9) $$H_{(-\infty,-1]} \cap H_{[0,\infty)} = \{0\}.$$

PROOF. First we assume (2.4). Then

$$H_{(-\infty,-1]} \cap H_{[0,\infty)} \subset H_{(-\infty,-1]} \cap H_{[-1,\infty)} = H_{\{-1\}},$$
$$H_{(-\infty,-1]} \cap H_{[0,\infty)} \subset H_{(-\infty,0]} \cap H_{[0,\infty)} = H_{\{0\}},$$

while $X_{-1}$ and $X_0$ are linearly independent since $\{X_n\}$ is nondeterministic. Thus, (2.9) follows.

Next we assume (2.9). By the arguments in [12], page 6, we see that $\{X_n\}$ is PND, whence, in particular, nondeterministic. Let $n \in \mathbf{N}$ and $X \in H_{(-\infty,-1]} \cap H_{[-n,\infty)}$. By the Regression Lemma, $X$ has a decomposition $X = Y + Z$ with $Y \in H_{(-\infty,-n-1]}$ and $Z \in H_{[-n,-1]}$. Then $Y = X - Z \in H_{(-\infty,-n-1]} \cap H_{[-n,\infty)}$. However, (2.9) implies $H_{(-\infty,-n-1]} \cap H_{[-n,\infty)} = \{0\}$, so that $Y = 0$ or $X = Z \in H_{[-n,-1]}$. Thus, $H_{(-\infty,-1]} \cap H_{[-n,\infty)} \subset H_{[-n,-1]}$. Since the converse implication $\supset$ is trivial, we obtain (1.2), whence (2.4). □

A stationary process $\{X_n\}$ is said to be *completely nondeterministic* if (2.9) holds. Thus, Theorem 2.3 asserts the equivalence between (2.4) and the complete nondeterminism of $\{X_n\}$. Complete nondeterminism was introduced by Sarason [25].

REMARK 2. In the first version of this manuscript, we raised the characterization of (2.4) in terms of the spectral density $\Delta(\cdot)$ as an open problem after remarking that (2.4) implies (2.9), whence that $\{X_n\}$ is PND. In the summer of 2004, Mohsen Pourahmadi, and then an anonymous referee, suggested the equivalence between (2.4) and (2.9), and both cited Bloomfield, Jewell and Hayashi [5], in which several characterizations of complete nondeterminism (2.9), in terms of the outer function determined by $\Delta(\cdot)$, which is essentially the same as $h(z)$ in (2.11) below, are given. Thus, we owe much of Theorem 2.3 to them.

2.2. *Representation in terms of MA and AR coefficients.* In this section we assume that the stationary process $\{X_n\}$ is purely nondeterministic. For $n \in \mathbf{N}$ and $m \in \mathbf{N} \cup \{0\}$, we can express the $(m+1)$-step predictor $P_{[-n,-1]}X_m$ uniquely in the form

$$(2.10) \qquad P_{[-n,-1]}X_m = \sum_{j=1}^{n} \phi_{n,j}^m X_{-j}.$$

We are concerned with representation of the real coefficients $\phi_{n,j}^m$, which we call the $(m+1)$-*step finite predictor coefficients*. In the 1-step case $m=0$, we have $\phi_{n,j}^0 = \phi_{n,j}$ by (1.1).

We consider the outer function

$$(2.11) \quad h(z) := \sqrt{2\pi} \exp\left\{\frac{1}{4\pi} \int_{-\pi}^{\pi} \frac{e^{i\lambda} + z}{e^{i\lambda} - z} \log \Delta(\lambda) \, d\lambda\right\}, \qquad z \in \mathbf{C}, \ |z| < 1.$$

The function $h(z)$ is holomorphic and has no zeros in $|z| < 1$, and it satisfies $2\pi\Delta(\lambda) = |h(e^{i\lambda})|^2$ a.e., where $h(e^{i\lambda}) := \lim_{r\uparrow 1} h(re^{i\lambda})$. We define the *MA*

*coefficients* $c_n$ by

$$h(z) = \sum_{n=0}^{\infty} c_n z^n, \qquad |z| < 1,$$

and the *AR coefficients* $a_n$ by

$$-1/h(z) = \sum_{n=0}^{\infty} a_n z^n, \qquad |z| < 1$$

(see Section 2 in [15]). Both $\{c_n\}$ and $\{a_n\}$ are real sequences, and we have $c_0 > 0$ and $\sum_0^\infty (c_n)^2 < \infty$. The coefficients $c_n$ and $a_n$ are actually those that appear in the following MA($\infty$) and AR($\infty$) representations, respectively, of $\{X_n\}$ [under suitable condition such as (2.15) below for the latter]:

$$(2.12) \qquad X_n = \sum_{j=-\infty}^{n} c_{n-j} \xi_j, \qquad n \in \mathbf{Z},$$

$$(2.13) \qquad \sum_{j=-\infty}^{n} a_{n-j} X_j + \xi_n = 0, \qquad n \in \mathbf{Z},$$

where $\{\xi_k\}$ is the innovation process given by $\xi_k = \varepsilon_k/\|\varepsilon_k\|$ with $\varepsilon_k$ in (2.8); see, for example, Chapter II in [23] for (2.12), and (4.9) in [15] for (2.13). By the assumption that $\{X_k\}$ is PND, $\{\xi_k\}$ forms a complete orthonormal system of $H$ such that, for every $n \in \mathbf{Z}$, the closed linear span of $\{\xi_k : -\infty < k \leq n\}$ in $H$ is equal to $H_{(-\infty, n]}$. Notice that the sums in (2.13) may not converge in norm in $H$.

EXAMPLE 2.4. Let $r \in (-1, 1)$. We consider the unique *causal* solution $X_n = \sum_{j=-\infty}^{n} r^{n-j} e_j$ to the AR(1) equation $X_n = r X_{n-1} + e_n$, where $\{e_n : n \in \mathbf{Z}\}$ is white noise, that is, a sequence in $H$ such that $(e_n, e_m) = \delta_{nm}$ (see, e.g., Section 4.1.1 in [22]). By standard computations, we find the equalities

$$\xi_n = e_n, \qquad \gamma(n) = \frac{r^{|n|}}{1 - r^2}, \qquad \Delta(\lambda) = \frac{1}{2\pi} \frac{1}{|1 - re^{i\lambda}|^2}, \qquad h(z) = \frac{1}{1 - rz},$$

$$c_n = r^n \quad (n \geq 0), \qquad a_0 = -1, \qquad a_1 = r, \qquad a_n = 0 \quad (n \geq 2).$$

We put

$$b_j^m := \sum_{k=0}^{m} c_k a_{j+m-k}, \qquad m, j \in \mathbf{N} \cup \{0\}.$$

In particular, $b_j^0 = c_0 a_j$. For $n \in \mathbf{N}$ and $m, j \in \mathbf{N} \cup \{0\}$, we define $b_k^m(n, j)$ recursively by

$$(2.14) \quad \begin{aligned} b_1^m(n, j) &= b_j^m, \\ b_{k+1}^m(n, j) &= \sum_{m_1=0}^{\infty} b_{n+1+m_1}^m b_k^{m_1}(n, j), \qquad k = 1, 2, \ldots. \end{aligned}$$

From the proof of Theorem 2.5 below, we see that, under the condition

$$\sum_{n=0}^{\infty} |a_n| < \infty, \tag{2.15}$$

which ensures the absolute convergence of the sums in (2.13), the sums in (2.14) also converge absolutely. We put, for $m \in \mathbf{N} \cup \{0\}$, $n \in \mathbf{N}$ and $j = 1, 2, \ldots, n$,

$$g_k^m(n,j) := \begin{cases} b_k^m(n,j), & k = 1, 3, \ldots, \\ b_k^m(n, n+1-j), & k = 2, 4, \ldots. \end{cases}$$

We write $\sum^{\infty-}$ for the improper sum: $\sum^{\infty-} = \lim_{M \to \infty} \sum^{M}$. The following theorem gives an explicit representation of the $(m+1)$-step finite predictor coefficients $\phi_{n,j}^m$ in (2.10), in terms of the MA and AR coefficients, under the absolute convergence of the sums in (2.13).

THEOREM 2.5. *We assume that the AR coefficients $a_n$ of a purely nondeterministic stationary process $\{X_n\}$ satisfy* (2.15). *Then we have $\phi_{n,j}^m = \sum_{k=1}^{\infty-} g_k^m(n,j)$ for $n \in \mathbf{N}$, $m \in \mathbf{N} \cup \{0\}$ and $j = 1, \ldots, n$, that is,*

$$P_{[-n,-1]} X_m = \sum_{j=1}^{n} \left\{ \sum_{k=1}^{\infty-} g_k^m(n,j) \right\} X_{-j}.$$

PROOF. For $m \in \mathbf{N} \cup \{0\}$ and $n \in \mathbf{N}$, we have the Wiener prediction formulas (see, e.g., Theorem 4.4 in [15])

$$P_{(-\infty,-1]} X_m = \sum_{j=1}^{\infty} b_j^m X_{-j}, \tag{2.16}$$

$$P_{[-n,\infty)} X_{-n-1-m} = \sum_{j=1}^{\infty} b_j^m X_{-n-1+j}, \tag{2.17}$$

the sums converging absolutely in $H$. Recall $P_n^k$ from (2.2). From (2.16), we have

$$P_n^1 X_m = \sum_{j=1}^{n} g_1^m(n,j) X_{-j} + \sum_{m_1=0}^{\infty} b_{n+1+m_1}^m X_{-n-1-m_1}.$$

From this and (2.17), it follows that

$$P_n^2 P_n^1 X_m = \sum_{j=1}^{n} g_1^m(n,j) X_{-j} + \sum_{m_1=0}^{\infty} b_{n+1+m_1}^m \sum_{j=1}^{\infty} b_j^{m_1} X_{-n-1+j}$$

$$= \sum_{j=1}^{n} \{g_1^m(n,j) + g_2^m(n,j)\} X_{-j}$$

$$+ \sum_{m_1=0}^{\infty} b_{n+1+m_1}^{m} \sum_{m_2=0}^{\infty} b_{n+1+m_2}^{m_1} X_{m_2}.$$

Similarly,

$$\begin{aligned}
P_n^3 P_n^2 P_n^1 X_m &= \sum_{j=1}^{n} \{g_1^m(n,j) + g_2^m(n,j)\} X_{-j} \\
&\quad + \sum_{m_1=0}^{\infty} b_{n+1+m_1}^{m} \sum_{m_2=0}^{\infty} b_{n+1+m_2}^{m_1} \sum_{j=1}^{\infty} b_j^{m_2} X_{-j} \\
&= \sum_{j=1}^{n} \{g_1^m(n,j) + g_2^m(n,j) + g_3^m(n,j)\} X_{-j} \\
&\quad + \sum_{m_1=0}^{\infty} b_{n+1+m_1}^{m} \sum_{m_2=0}^{\infty} b_{n+1+m_2}^{m_1} \sum_{m_3=0}^{\infty} b_{n+1+m_3}^{m_2} X_{-n-1-m_3}.
\end{aligned}$$

Repeating this argument, we see that $\phi_{n,j}^k(X_m)$ in Lemma 2.1 with $Y = X_m$ are given by $\phi_{n,j}^k(X_m) = \sum_{l=1}^{k} g_l^m(n,j)$. The condition (2.15) implies $\sum_{0}^{\infty}(a_n)^2 < \infty$, whence (1.3) (see, e.g., Proposition 4.2 in [15]). Thus, the theorem follows from Theorem 2.2. □

2.3. *Representation by absolutely convergent series.* In the applications which we discuss later, the finite predictor coefficients $\phi_{n,j}$ in (1.1) need to be expressed by an absolutely convergent series made up of $a_k$ and $c_k$. In this section we first give such an expression for $b_k^m(n,j)$. In the 1-step case $m = 0$, the result yields the desired representation for $\phi_{n,j}$.

We write $\mathcal{R}_0$ for the class of *slowly varying functions* at infinity: the class of positive, measurable $\ell(\cdot)$, defined on some neighborhood $[A, \infty)$ of infinity, such that $\lim_{x \to \infty} \ell(\lambda x)/\ell(x) = 1$ for all $\lambda > 0$ (see Chapter 1 in [4] for background).

Throughout this section we assume that the stationary process $\{X_n\}$ satisfies one of the following conditions (A1) and (A2):

(A1) $\{X_n\}$ is purely nondeterministic, and $\{a_n\}$ and $\{c_n\}$ satisfy, respectively, (2.15) and

$$\sum_{n=0}^{\infty} |c_n| < \infty. \tag{2.18}$$

(A2) $\{X_n\}$ is purely nondeterministic and, for $d \in (0, 1/2)$ and $\ell(\cdot) \in \mathcal{R}_0$, $\{c_n\}$ and $\{a_n\}$ satisfy, respectively,

$$c_n \sim n^{-(1-d)} \ell(n), \qquad n \to \infty, \tag{2.19}$$

$$a_n \sim n^{-(1+d)} \frac{1}{\ell(n)} \cdot \frac{d \sin(\pi d)}{\pi}, \qquad n \to \infty. \tag{2.20}$$

It should be noticed that (2.20) implies (2.15).

In this paper we say that a stationary process $\{X_n\}$ has *long memory* (resp. *short memory*) if $\sum_{k=-\infty}^{\infty} |\gamma(k)| = \infty$ (resp. $< \infty$). See [2], page 6, and Section 13.2 in [6]. By (2.12), the autocovariance function $\gamma(\cdot)$ has the expression

$$\gamma(n) = \sum_{k=0}^{\infty} c_{|n|+k} c_k, \qquad n \in \mathbf{Z}. \tag{2.21}$$

Hence, (2.18) implies that

$$\sum_{n=0}^{\infty} |\gamma(n)| \leq \left(\sum_{k=0}^{\infty} |c_k|\right)^2 < \infty.$$

Thus, $\{X_n\}$ has short memory under (A1). On the other hand, by (2.21) and [14], Proposition 4.3, (2.19) implies that

$$\gamma(n) \sim n^{-(1-2d)} \ell(n)^2 B(d, 1-2d), \qquad n \to \infty. \tag{2.22}$$

Since $0 < 1 - 2d < 1$, we see that $\{X_n\}$ has long memory under (A2). We remark that, under suitable conditions, (2.19), (2.20) and (2.22) are equivalent (see Theorem 5.1 in [15]).

EXAMPLE 2.6. For $d \in (-1/2, 1/2)$ and $p, q \in \mathbf{N} \cup \{0\}$, a stationary process $\{X_n\}$ is said to be a fractional ARIMA$(p, d, q)$ process if it has a spectral density $\Delta(\cdot)$ of the form

$$\Delta(\lambda) = \frac{1}{2\pi} \frac{|\theta(e^{i\lambda})|^2}{|\phi(e^{i\lambda})|^2} |1 - e^{i\lambda}|^{-2d}, \qquad -\pi \leq \lambda \leq \pi,$$

where $\phi(z)$ and $\theta(z)$ are polynomials with real coefficients of degrees $p$ and $q$, respectively. We assume that $\phi(z)$ and $\theta(z)$ have no common zeros, and that neither $\phi(z)$ nor $\theta(z)$ has zeros in the closed unit disk $\{z \in \mathbf{C} : |z| \leq 1\}$. We also assume without loss of generality that $\theta(0)/\phi(0) > 0$. Then the outer function $h(\cdot)$ is given by $h(z) = (1-z)^{-d} \theta(z)/\phi(z)$ (see, e.g., Section 2 in [17]). If $0 < d < 1/2$, then $\{X_n\}$ satisfies (A2) for some constant function $\ell(\cdot)$ (see Corollary 3.1 in [19]). If $d = 0$, then $\{X_n\}$ is also called an ARMA$(p, q)$ process (see Chapter 3 in [6]), and both $\{c_n\}$ and $\{a_n\}$ decay exponentially, whence (A1) is satisfied.

We put

$$B_n := \sum_{v=0}^{\infty} |c_v a_{n+v}|, \qquad n \in \mathbf{N} \cup \{0\}.$$

For $n, k, u, v \in \mathbf{N} \cup \{0\}$, we define $D_k(n, u, v)$ recursively by

$$D_0(n, u, v) := \delta_{uv},$$

$$D_{k+1}(n, u, v) := \sum_{w=0}^{\infty} B_{n+v+w} D_k(n, u, w).$$

We have, for example,

$$D_3(n, u, v) = \sum_{v_1=0}^{\infty} \sum_{v_2=0}^{\infty} B_{n+v+v_1} B_{n+v_1+v_2} B_{n+v_2+u}.$$

By the Fubini–Tonelli theorem, we have $D_k(n, u, v) = D_k(n, v, u)$.

LEMMA 2.7. *We assume either* (A1) *or* (A2). *Then, for* $k, n, v \in \mathbf{N} \cup \{0\}$,

$$\sum_{u=0}^{\infty} D_k(n, u, v) < \infty \quad \text{and} \quad \sum_{u=0}^{\infty} D_k(n, u, v)^2 < \infty,$$

*respectively. In particular, we have* $D_k(n, u, v) < \infty$ *for* $k, n, u, v \in \mathbf{N} \cup \{0\}$.

PROOF. First we assume (A1). Then

$$\sum_{m=0}^{\infty} B_m \leq \left\{ \sum_{u=0}^{\infty} |c_u| \right\} \left\{ \sum_{u=0}^{\infty} |a_u| \right\} < \infty.$$

This and the nonnegativity of $B_m$ imply, for example,

$$\sum_{u=0}^{\infty} D_3(n, u, v) = \sum_{u=0}^{\infty} \sum_{v_1=0}^{\infty} \sum_{v_2=0}^{\infty} B_{n+v+v_1} B_{n+v_1+v_2} B_{n+v_2+u}$$

$$\leq \left\{ \sum_{m=0}^{\infty} B_m \right\}^3 < \infty.$$

The general case can be proved in the same way.

Next we assume (A2). The proof in this case is the same as that of Lemma 2.1 in [16]. By (A2) and Proposition 4.3 in [14], we have $B_n = O(n^{-1})$ as $n \to \infty$. Therefore, for $n \in \mathbf{N}$, $f_u \mapsto \sum_{v=0}^{\infty} B_{n+u+v} f_v$ defines a bounded linear operator on $l^2$ (see Chapter IX in [11]). Since $D_{k+1}(n, u, v) = \sum_w B_{n+u+w} D_k(n, w, v)$, we obtain the desired result by induction on $k$. □

We put

(2.23) $$\beta_n := \sum_{v=0}^{\infty} c_v a_{v+n}, \qquad n = 0, 1, \ldots.$$

In view of Lemma 2.7, we may define $\delta_k(n,u,v)$ recursively by, for $k,n,u,v \in \mathbf{N} \cup \{0\}$,

(2.24)
$$\delta_0(n,u,v) = \delta_{uv},$$

$$\delta_{k+1}(n,u,v) = \sum_{w=0}^{\infty} \beta_{n+v+w} \delta_k(n,u,w).$$

By Lemma 2.7 and the Fubini theorem, we have $\delta_k(n,u,v) = \delta_k(n,v,u)$.

The following theorem expresses $b_k^m(n,j)$ as an absolutely convergent series.

THEOREM 2.8. *We assume either* (A1) *or* (A2). *Then, for* $n,k \in \mathbf{N}$ *and* $m,j \in \mathbf{N} \cup \{0\}$,

(2.25)
$$b_k^m(n,j) = \sum_{v=0}^{m} c_{m-v} \sum_{u=0}^{\infty} a_{j+u} \delta_{k-1}(n+1,u,v),$$

*the sum converging absolutely.*

PROOF. By Lemma 2.7 and (2.15), we have

(2.26)
$$\sum_{u=0}^{\infty} |a_{j+u}| D_{k-1}(n+1,u,v)$$
$$\leq \left\{\sup_u D_{k-1}(n+1,u,v)\right\} \sum_{u=0}^{\infty} |a_{j+u}| < \infty.$$

Thus, the right-hand side of (2.25), which we denote by $B_k^m(n,j)$, converges absolutely. To prove the proposition, it is enough show that $B_k^m(n,j)$ satisfies the same recursion as (2.14).

First we have

$$B_1^m(n,j) = \sum_{v=0}^{m} c_{m-v} \sum_{u=0}^{\infty} a_{j+u} \delta_{uv} = \sum_{v=0}^{m} c_{m-v} a_{j+v} = b_j^m,$$

as desired. Next, the Fubini–Tonelli theorem and (2.26) yield, for $k \geq 1$,

$$\sum_{u=0}^{\infty} a_{j+u} \delta_k(n+1,u,v)$$
$$= \sum_{u=0}^{\infty} a_{j+u} \sum_{w=0}^{\infty} \left\{\sum_{m_1=w}^{\infty} c_{m_1-w} a_{n+1+v+m_1}\right\} \delta_{k-1}(n+1,u,w)$$
$$= \sum_{m_1=0}^{\infty} a_{n+1+v+m_1} \sum_{w=0}^{m_1} c_{m_1-w} \sum_{u=0}^{\infty} a_{j+u} \delta_{k-1}(n+1,u,w)$$

$$= \sum_{m_1=0}^{\infty} a_{n+1+v+m_1} B_k^{m_1}(n,j),$$

so that

$$B_{k+1}^m(n,j) = \sum_{v=0}^{m} c_{m-v} \sum_{m_1=0}^{\infty} a_{n+1+v+m_1} B_k^{m_1}(n,j)$$

$$= \sum_{m_1=0}^{\infty} \left\{ \sum_{v=0}^{m} c_{m-v} a_{n+1+m_1+v} \right\} B_k^{m_1}(n,j)$$

$$= \sum_{m_1=0}^{\infty} b_{n+1+m_1}^m B_k^{m_1}(n,j).$$

Thus, $B_k^m(n,j)$ satisfies (2.14). $\square$

For applications in later sections, we consider the case $m=0$ separately. We put

$$d_k(n,j) := \delta_k(n,0,j), \qquad n,k,j \in \mathbf{N} \cup \{0\}.$$

Then, by (2.24), $d_k(n,j)$ satisfies the following recursion: for $k,n,j \in \mathbf{N} \cup \{0\}$,

$$d_0(n,j) = \delta_{j0},$$

(2.27)

$$d_{k+1}(n,j) = \sum_{v=0}^{\infty} \beta_{n+j+v} d_k(n,v).$$

More explicitly, $d_k(n,j)$ are given by, for $n,j \in \mathbf{N} \cup \{0\}$,

$$d_1(n,j) = \beta_{n+j}, \qquad d_2(n,j) = \sum_{v_1=0}^{\infty} \beta_{n+j+v_1} \beta_{n+v_1},$$

and, for $k=3,4,\ldots$,

$$d_k(n,j) = \sum_{v_1=0}^{\infty} \cdots \sum_{v_{k-1}=0}^{\infty} \beta_{n+j+v_{k-1}} \beta_{n+v_{k-1}+v_{k-2}} \cdots \beta_{n+v_2+v_1} \beta_{n+v_1},$$

the sums converging absolutely.

We put

$$b_k(n,j) := b_k^0(n,j), \qquad g_k(n,j) := g_k^0(n,j)$$

for $(k,n,j)$, for which the right-hand sides are defined. Then, for $n \in \mathbf{N}$ and $j = 1,2,\ldots,n$, we have

(2.28) $$g_k(n,j) = \begin{cases} b_k(n,j), & k=1,3,\ldots, \\ b_k(n,n+1-j), & k=2,4,\ldots. \end{cases}$$

By Theorems 2.5 and 2.8, we immediately obtain the following final form of the representation of the 1-step finite predictor coefficients $\phi_{n,j}$.

THEOREM 2.9. *We assume either* (A1) *or* (A2). *Then, for $n \in \mathbf{N}$ and $j = 1, \ldots, n$, we have $\phi_{n,j} = \sum_{k=1}^{\infty-} g_k(n,j)$ with (2.28) and*

$$b_1(n,v) = c_0 a_v, \qquad v \geq 0,$$

$$b_k(n,v) = c_0 \sum_{u=0}^{\infty} a_{v+u} d_{k-1}(n+1, u), \qquad k \geq 2, \ v \geq 0,$$

*the sum on the right-hand side converging absolutely.*

**3. The rate of convergence of finite predictor coefficients.** If the stationary process $\{X_n\}$ is PND and satisfies (2.15), then we have the Wiener prediction formula (2.16) with $m = 0$ or (1.4) with

(3.1) $$\phi_j = c_0 a_j, \qquad j \in \mathbf{N}.$$

We call $\phi_j$ the *infinite predictor coefficients*. It holds that

$$\lim_{n \to \infty} \phi_{n,j} = \phi_j, \qquad j \in \mathbf{N}$$

(see, e.g., Theorem 7.14 in [22]). In this section we investigate the rate for long memory processes at which $\phi_{n,j}$ converges to $\phi_j$. Notice that, by (2.14), (2.28) and (3.1), we have

(3.2) $$\phi_j = b_1(n,j) = g_1(n,j), \qquad n \in \mathbf{N}, \ j = 1, \ldots, n.$$

Thus, $\phi_j$ is the first term of the series $\sum_{k=1}^{\infty-} g_k(n,j)$ in Theorem 2.9 expressing $\phi_{n,j}$. This suggests the usefulness of the expression for our purpose.

Throughout this section, we assume that the stationary process $\{X_n\}$ satisfies (A2) in Section 2.3 (long memory).

For $u \geq 0$, we put

$$f_1(u) := \frac{1}{\pi(1+u)}, \qquad f_2(u) := \frac{1}{\pi^2} \int_0^\infty \frac{ds_1}{(s_1+1)(s_1+1+u)},$$

and, for $k = 3, 4, \ldots$,

$$f_k(u) := \frac{1}{\pi^k} \int_0^\infty ds_{k-1} \cdots \int_0^\infty ds_1 \frac{1}{(s_{k-1}+1)}$$

$$\times \left\{ \prod_{m=1}^{k-2} \frac{1}{(s_{m+1}+s_m+1)} \right\} \frac{1}{(s_1+1+u)}$$

(see Section 3 in [17]; see also Section 6 in [15]).

LEMMA 3.1. (i) $\sum_{k=1}^{\infty} f_{2k}(0) x^{2k} = (\pi^{-1} \arcsin x)^2$ *for* $|x| < 1$;
(ii) $\sum_{k=1}^{\infty} f_{2k-1}(0) x^{2k-1} = \pi^{-1} \arcsin x$ *for* $|x| < 1$.

PROOF. Let $j \geq 1$. We easily see that $f_{1+j}(u) = \int_0^\infty f_1(s+u) f_j(s)\, ds$ for $u \geq 0$. Hence, we have, for $u \geq 0$,

$$\begin{aligned} f_{2+j}(0) &= \int_0^\infty f_1(s) f_{j+1}(s)\, ds \\ &= \int_0^\infty ds\, f_1(s) \int_0^\infty f_1(u+s) f_j(u)\, du \\ &= \int_0^\infty \left\{ \int_0^\infty f_1(s+u) f_1(s)\, ds \right\} f_j(u)\, du \\ &= \int_0^\infty f_2(u) f_j(u)\, du. \end{aligned}$$

Repeating this argument, we obtain

$$(3.3) \qquad \int_0^\infty f_i(u) f_j(u)\, du = f_{i+j}(0), \qquad i, j \in \mathbf{N}.$$

Thus, the assertion (i) follows from Lemma 6.5 in [15], while (ii) follows from Lemma 3.4 in [16]. □

Recall $d_k(n, u)$ from Section 2.3.

PROPOSITION 3.2. (i) *For $r \in (1, \infty)$, there exists $N \in \mathbf{N}$ such that*

$$(3.4)\ 0 < d_k(n, u) \leq \frac{f_k(0)\{r \sin(\pi d)\}^k}{n}, \qquad u \in \mathbf{N} \cup \{0\},\ k \in \mathbf{N},\ n \geq N.$$

(ii) *For $k \in \mathbf{N}$ and $u \in \mathbf{N} \cup \{0\}$, $d_k(n, u) \sim n^{-1} f_k(0) \sin^k(\pi d)$ as $n \to \infty$.*

PROOF. Let $r > 1$. Recall $\beta_n$ from (2.23). The condition (A2) implies

$$(3.5) \qquad \beta_n \sim \frac{\sin(\pi d)}{\pi} n^{-1}, \qquad n \to \infty$$

(see Proposition 4.3 in [14]). Thus, for $n$ large enough,

$$0 < \beta_{[ns]+n+u} \leq \frac{r^{1/2} \sin(\pi d)}{\pi([ns] + n + u)}, \qquad s \geq 0,\ u \in \mathbf{N} \cup \{0\}.$$

Since we have, for $n$ large enough,

$$\frac{1}{[ns] + n + u} \leq \frac{r^{1/2}}{n(s+1)}, \qquad s \geq 0,\ u \in \mathbf{N} \cup \{0\},$$

there exists $N_1 \in \mathbf{N}$ such that

$$(3.6) \quad 0 < \beta_{[ns]+n+u} \leq \frac{r \sin(\pi d)}{\pi(s+1)} n^{-1}, \qquad s \geq 0,\ u \in \mathbf{N} \cup \{0\},\ n \geq N_1.$$

In the same way, we can choose $N_2$ so that

$$(3.7) \quad 0 < \beta_{[ns_2]+[ns_1]+n} \leq \frac{r\sin(\pi d)}{\pi(s_2+s_1+1)} n^{-1}, \qquad s_1, s_2 \geq 0, \ n \geq N_2.$$

Therefore, we have, for $n \geq N := \max(N_1, N_2)$,

$$\begin{aligned}
0 < d_3(n,u) &= \int_0^\infty ds_2 \int_0^\infty ds_1 \cdot \beta_{[s_2]+n} \cdot \beta_{[s_2]+[s_1]+n} \cdot \beta_{[s_1]+n+u} \\
&= n^2 \int_0^\infty ds_2 \int_0^\infty ds_1 \cdot \beta_{[ns_2]+n} \cdot \beta_{[ns_2]+[ns_1]+n} \cdot \beta_{[ns_1]+n+u} \\
&\leq \frac{\{r\sin(\pi d)\}^3}{n} \frac{1}{\pi^3} \int_0^\infty ds_2 \int_0^\infty ds_1 \frac{1}{(s_2+1)(s_2+s_1+1)(s_1+1)} \\
&= \frac{\{r\sin(\pi d)\}^3}{n} f_3(0),
\end{aligned}$$

which implies (3.4) with $k = 3$. Notice that $N$ is independent of the choice $k = 3$. We can prove (3.4) for general $k$ and the same $N$ in a similar fashion.

We also prove (ii) only for $k = 3$; the general case can be treated in the same way. By (3.5), we have

$$(3.8) \qquad \lim_{n\to\infty} n\beta_{[ns]+n+u} = \frac{\sin(\pi d)}{\pi(s+1)}, \qquad s \geq 0, \quad u \in \mathbf{N} \cup \{0\},$$

$$(3.9) \quad \lim_{n\to\infty} n\beta_{[ns_2]+[ns_1]+n} = \frac{\sin(\pi d)}{\pi(s_2+s_1+1)}, \qquad s_1, s_2 \geq 0.$$

By (3.6)–(3.9) and the dominated convergence theorem, we obtain

$$\begin{aligned}
\lim_{n\to\infty} \int_0^\infty ds_2 \int_0^\infty ds_1 \cdot n\beta_{[ns_2]+n} \cdot n\beta_{[ns_2]+[ns_1]+n} \cdot n\beta_{[ns_1]+n+u} \\
= \frac{\sin^3(\pi d)}{\pi^3} \int_0^\infty ds_2 \int_0^\infty ds_1 \frac{1}{(s_2+1)(s_2+s_1+1)(s_1+1)}.
\end{aligned}$$

This implies $\lim_n nd_3(n,u) = \sin^3(\pi d) f_3(0)$ or (ii) with $k = 3$, as desired. □

The following theorem gives the rate for long memory processes at which $\phi_{n,j}$ converges to $\phi_j$. It applies, in particular, to the fractional ARIMA$(p, d, q)$ processes with $0 < d < 1/2$.

THEOREM 3.3. *We assume* (A2). *Then we have, for $j \in \mathbf{N}$,*

$$\lim_{n\to\infty} n\{\phi_{n,j} - \phi_j\} = d^2 \sum_{u=j}^\infty \phi_u.$$

PROOF. Let $r > 1$ be chosen so that $0 < r\sin(\pi d) < 1$. By Lemma 3.1,

$$\text{(3.10)} \qquad \sum_{k=1}^{\infty} f_k(0)\{r\sin(\pi d)\}^k < \infty.$$

Let $N$ be such that (3.4) holds. Then, for $n \geq N$ and $j = 1, \ldots, n$,

$$n\left|\sum_{u=0}^{\infty} a_{n-j+u} \sum_{k=1}^{\infty} d_{2k-1}(n,u)\right|$$

$$\leq \sum_{u=0}^{\infty} |a_{n-j+u}| \sum_{k=1}^{\infty} n d_{2k-1}(n,u)$$

$$\leq \left[\sum_{k=1}^{\infty} f_{2k-1}(0)\{r\sin(\pi d)\}^{2k-1}\right] \sum_{u=n-j}^{\infty} |a_u|,$$

so that

$$\lim_{n\to\infty} n \sum_{u=0}^{\infty} a_{n-j+u} \sum_{k=1}^{\infty} d_{2k-1}(n,u) = 0.$$

Proposition 3.2, Lemma 3.1 and the dominated convergence theorem yield

$$\lim_{n\to\infty} n \sum_{u=0}^{\infty} a_{j+u} \sum_{k=1}^{\infty} d_{2k}(n,u)$$

$$= \left\{\sum_{k=1}^{\infty} f_{2k}(0) \sin^{2k}(\pi d)\right\} \sum_{u=0}^{\infty} a_{j+u} = d^2 \sum_{u=j}^{\infty} a_u.$$

Therefore, by Theorem 2.9 and (3.2) we have

$$\lim_{n\to\infty} n\{\phi_{n-1,j} - \phi_j\}$$

$$= \lim_{n\to\infty} \left\{n \sum_{k=1}^{\infty} b_{2k+1}(n-1,j) + n \sum_{k=1}^{\infty} b_{2k}(n-1, n-j)\right\}$$

$$= \lim_{n\to\infty} c_0 n \sum_{u=0}^{\infty} a_{j+u} \sum_{k=1}^{\infty} d_{2k}(n,u) + \lim_{n\to\infty} c_0 n \sum_{u=0}^{\infty} a_{n-j+u} \sum_{k=1}^{\infty} d_{2k-1}(n,u)$$

$$= c_0 d^2 \sum_{u=j}^{\infty} a_u = d^2 \sum_{u=j}^{\infty} \phi_u,$$

as desired.  □

The next proposition, which follows from the proof of Theorem 3.3, shows that, under (A2), the sum $\sum_{k=1}^{\infty} g_k(n,j)$ converges absolutely for $n$ large enough and $j = 1, \ldots, n$.

PROPOSITION 3.4. *We assume* (A2). *If $N$ is such that* (3.4) *holds for some $r \in (1, \infty)$, then $\sum_{k=1}^{\infty} |g_k(n,j)| < \infty$ for $n \geq N$ and $j = 1, \ldots, n$.*

For processes with short memory, we have the following result.

PROPOSITION 3.5. *We assume* (A1). *If $N$ is such that $(\sum_{j=0}^{\infty} |c_j|) \times (\sum_{k=N+1}^{\infty} |a_k|) < 1$, then $\sum_{k=1}^{\infty} |g_k(n,j)| < \infty$ for $n \geq N$ and $j = 1, \ldots, n$.*

We omit the proof.

**4. Baxter's inequality for long memory processes.** In this section we prove Baxter's inequality (1.5).

THEOREM 4.1. *We assume* (A2). *Then there exists a positive constant $M$ such that* (1.5) *holds for all $n \in \mathbf{N}$.*

PROOF. Let $r > 1$ be chosen so that $0 < r \sin(\pi d) < 1$. Then we have (3.10). By Proposition 3.2 and (2.20), we may take a positive integer $N$ such that both (3.4) and $a_n > 0$ hold for $n \geq N$. Pick $\delta \in (0, d)$. By (2.20) and [4], Theorem 1.5.6(iii) (Potter-type bounds), we may assume that

$$(4.1) \quad a_m/a_n \leq 2\max\{(n/m)^{1+d-\delta}, (m/n)^{1+d+\delta}\}, \qquad m, n \geq N.$$

By Theorem 2.9 and (3.2), we have, for $n \geq N + 3$,

$$\sum_{j=1}^{n-1} |\phi_{n-1,j} - \phi_j| \leq c_0 \sum_{j=1}^{n-1} \sum_{u=0}^{\infty} |a_{u+j}| \sum_{k=1}^{\infty} d_{2k}(n,u)$$

$$+ c_0 \sum_{j=1}^{n-1} \sum_{u=0}^{\infty} |a_{u+n-j}| \sum_{k=1}^{\infty} d_{2k-1}(n,u)$$

$$= c_0 \sum_{k=1}^{\infty} \sum_{u=0}^{\infty} d_k(n,u) \sum_{j=1}^{n-1} |a_{u+j}|$$

$$= c_0 \{G_1(n) + G_2(n)\},$$

where

$$G_1(n) := \sum_{k=1}^{\infty} \sum_{u=0}^{\infty} d_k(n,u) \sum_{j=1}^{N+1} |a_{u+j}|,$$

$$G_2(n) := \sum_{k=1}^{\infty} \sum_{u=0}^{\infty} d_k(n,u) \sum_{j=N+2}^{n-1} a_{u+j}.$$

For $n \geq N+3$, we have

$$G_1(n) \leq n^{-1}\left[\sum_{k=1}^{\infty} f_k(0)\{r\sin(\pi d)\}^k\right] \sum_{j=1}^{N+1} \sum_{u=0}^{\infty} |a_{u+j}|,$$

and

$$G_2(n) \leq \sum_{k=1}^{\infty} \int_0^{\infty} du \cdot d_k(n,[u]) \int_N^n a_{[u]+[s]+2}\,ds$$

$$= na_n \sum_{k=1}^{\infty} \int_0^{\infty} du \cdot nd_k(n,[nu]) \int_{N/n}^1 \frac{a_{[nu]+[ns]+2}}{a_n}\,ds.$$

By (4.1), we have, for $u > 0$, $n \geq N+3$ and $N/n \leq s \leq 1$,

$$\frac{a_{[nu]+[ns]+2}}{a_n} \leq 2\max\left\{\left(\frac{n}{[nu]+[ns]+2}\right)^{1+\delta+d}, \left(\frac{n}{[nu]+[ns]+2}\right)^{1-\delta+d}\right\}$$

$$\leq 2\max\{(u+s)^{-(1+d+\delta)}, (u+s)^{-(1+d-\delta)}\}.$$

Hence, by (3.4),

$$G_2(n) \leq 2\int_0^1 ds \int_0^{\infty} \{(u+s)^{-(1+d+\delta)} + (u+s)^{-(1+d-\delta)}\}\,du$$

$$\times na_n\left[\sum_{k=1}^{\infty} f_k(0)\{r\sin(\pi d)\}^k\right]$$

$$\leq 2\left\{\frac{1}{(\delta+d)(1-d-\delta)} + \frac{1}{(d-\delta)(1-d+\delta)}\right\}$$

$$\times na_n\left[\sum_{k=1}^{\infty} f_k(0)\{r\sin(\pi d)\}^k\right].$$

Combining these estimates, we obtain

$$\limsup_{n\to\infty}\left\{n^d\ell(n)\sum_{j=1}^{n-1}|\phi_{n-1,j} - \phi_j|\right\} < \infty.$$

Since $\sum_{k=n}^{\infty}\phi_k = c_0\sum_{k=n}^{\infty}a_k \sim c_0\sin(\pi d)/\{\pi n^d\ell(n)\}$ as $n\to\infty$, the theorem follows. $\square$

**Acknowledgments.** We express our gratitude to Mohsen Pourahmadi and an anonymous referee for leading us to Theorem 2.3. We are also grateful to another referee and an Associate Editor whose constructive comments resulted in a restructuring of the paper.

DEPARTMENT OF MATHEMATICS
FACULTY OF SCIENCE
HOKKAIDO UNIVERSITY
SAPPORO 060-0810
JAPAN
E-MAIL: inoue@math.sci.hokudai.ac.jp
y-kasa@math.sci.hokudai.ac.jp